\documentclass[a4paper,12pt]{article}
\usepackage[utf8]{inputenc}
\title{On a connection between layer-adapted exponentially graded and S-type meshes}
\author{Sebastian Franz \footnote{%
                          Technische Universit\"{a}t Dresden, 
                          Institut f\"{u}r Numerische Mathematik,
                          01062 Dresden, Germany,
                          \mbox{e-mail}: sebastian.franz@tu-dresden.de}        
                        \footnote{
                          BTU Cottbus,
                          Institut f\"ur Mathematik, 
                          03046 Cottbus, Germany}
        \and
        Christos Xenophontos\footnote{%
                              Department of Mathematics and Statistics, University of Cyprus, P.O. Box
                              20537, 1678 Nicosia, Cyprus}}
\date{\today}

\usepackage{fancyhdr} 
\usepackage{enumerate}

\usepackage{amsmath}
\usepackage{amssymb}
\usepackage{amsthm}

\makeatletter
\let\my@saved@original@eqref\eqref 
\renewcommand*{\eqref}[1]{
  \begingroup
    \let\normalfont\relax
    \my@saved@original@eqref{#1}
  \endgroup
}
\makeatother

\usepackage{cite}

\usepackage{color}

\newcommand{\e}{\mathrm{e}}

\newcommand{\grad}{\nabla}

\newcommand{\eps}{\varepsilon}
\newcommand{\norm}[2]{\|{#1}\|_{#2}}

\newcommand{\tnorm}[1]{\left|\!\!\;\left|\!\!\;\left| {#1}
                       \right|\!\!\;\right|\!\!\;\right|}

\newcommand{\R}{\mathbb{R}}

\newcounter{tmp}
\newcommand{\makeballnumber}[1]{\setcounter{tmp}{\theenumi}%
\setcounter{enumi}{#1}%
\leavevmode \csname beamer@@tmpl@enumerate item\endcsname%
\setcounter{enumi}{\thetmp}}
\newcommand{\makeball}{\leavevmode \csname beamer@@tmpl@itemize item\endcsname}

\definecolor{seb}{rgb}{0.9,0,0}

\makeatletter
\renewcommand*\env@matrix[1][r]{\hskip -\arraycolsep
  \let\@ifnextchar\new@ifnextchar
  \array{*\c@MaxMatrixCols #1}}
\makeatother

\numberwithin{equation}{section}

\DeclareMathAlphabet{\mathcal}{OMS}{cmsy}{m}{n}

\theoremstyle{plain}

\begin{document}
 \pagestyle{fancy}
  \maketitle
  \begin{abstract}
    In this short note we analyse a connection between the exponentially graded and 
    the class of S-type meshes for singularly perturbed problems. As a by-product we 
    obtain a slightly modified and more general class of layer-adapted meshes.
  \end{abstract}



  \section{Introduction}
    In the numerical analysis of singularly perturbed problems the layer behaviour
    of the solution often creates challenging problems. One approach to deal theoretically
    and practically with these problems lies within using layer-adapted meshes.
    
    The idea is pretty old and goes back to Bakhvalov \cite{Bak69} for exponentially fitted 
    and Shishkin \cite{Shi90} for piecewise equidistant meshes. In \cite{RL99} a uniform 
    characterisation of meshes related to Shishkin meshes was introduced with the class of
    S-type meshes. Besides these meshes many more exist and numerical methods behave slightly 
    different on each of them, see \cite{RST08, Linss10}.
    
    Using $N$ as the number of cells in one space dimension we find in the literature
    convergence results for many methods on the previously mentioned meshes and singularly 
    perturbed problems, let us just mention the book \cite{RST08}. On S-type meshes they have 
    usually the form
    \[
      \norm{u-u_h}{}\leq C (N^{-1}\max|\psi'|)^k
    \]
    for some norm $\norm{\cdot}{}$, order $k$ and \emph{mesh characterisation function} $\psi$,
    where the constant $C$ here and further on is independent of $\eps$ and $N$.
    Meshes facilitating the estimate with $\max|\psi'|\leq C$ are called \emph{optimal meshes}, see 
    \cite{RTU14}. There are some meshes known with this desirable property, just to name the 
    Bakhvalov-S-mesh \cite{RL99} or those considered in \cite{RTU14}.
    
    Using a different fitting approach one finds the exponentially graded eXp-mesh \cite{Xenophontos02, CX14}.
    This is also an optimal mesh, see e.g. \cite{FrLX16}, but not an S-type mesh. In this note
    we will show a connection to the class of S-type meshes and thereby casting the eXp-mesh
    as an generalised S-mesh.

  \section{The eXp-mesh as a modified S-type mesh}
    We follow \cite{Xenophontos02, CX14} and define the eXp-mesh in one space dimension for a 
    convection-diffusion problem.
    Let $\sigma>0$ be a given parameter that depends on the polynomial degree of the underlying 
    discrete space. Then we define the constant
    \[
        C_{\sigma,\varepsilon}=1-\exp \left( -\frac{\beta }{\sigma\varepsilon }\right).
    \]
    Using this constant we obtain the mesh generating function
    \[
        \varphi_{eXp}(t)=-\ln \left[ 1-2C_{\sigma,\varepsilon }t\right] ,t\,\in \lbrack 0,1/2-1/N],
    \]
    and the transition point 
    \[
        \lambda_{eXp}=\frac{\sigma\eps}{\beta}\varphi_{eXp}\left( \frac{1}{2}-\frac{1}{N}\right)\leq\frac{1}{2},
    \]
    which includes an implicit assumption 
    $
        \eps\leq \frac{\beta}{2\sigma}(\ln (N/2))^{-1}.
    $
    Then the points of the mesh are given by
    \begin{gather}\label{eq:eXp}
        x_k = \begin{cases}
                  \frac{\sigma\eps}{\beta}\varphi_{eXp}\left(\frac{i}{N}\right),
                                              &i=0,\dots,N/2-1,\\
                  1-2(1-\lambda_{eXp})\frac{N-i}{N+2},&i=N/2-1,\dots,N.
              \end{cases}
    \end{gather}
    For an S-type mesh with mesh generating function $\phi:[0,1/2]\to\R_+$ the definitions look very similar. 
    Here it holds
    \[
       \lambda=\frac{\sigma\eps}{\beta}\ln N\leq\frac{1}{2}
    \]
    with the implicit assumption 
    $
        \eps\leq \frac{\beta}{2\sigma}(\ln N)^{-1}
    $
    and
    \begin{gather}\label{eq:Stype}
        x_k = \begin{cases}
                  \frac{\sigma\eps}{\beta}\phi\left(\frac{i}{N}\right),
                                              &i=0,\dots,N/2,\\
                  1-2(1-\lambda)\frac{N-i}{N},&i=N/2,\dots,N.
              \end{cases}
    \end{gather}
    Thus the difference lies in the mesh-generating function and the number of mesh-nodes in the layer region.
    
    For practical computations on a computer with machine precision $eps$, we obtain for small $\eps$,
    or more precise for
    \[
       \eps<\frac{\beta}{\sigma|\ln eps|},
    \]
    that $C_{\sigma,\eps}=1$. For example, having
    $
       \beta=1,\,\sigma=2,\,eps=2\cdot 10^{-16}
    $
    this happens for 
    $
       \eps\leq 0.0138.
    $
    Then the mesh generating function of the eXp-mesh simplifies to
    \[
       \tilde \varphi_{eXp}(t)=-\ln \left[ 1-2t\right] ,\,t\in \lbrack 0,1/2-1/N],
    \]
    or mapped onto the interval $[0,1/2]$ (in order to compare with S-type meshes)
    \[
       \phi_{eXp}(t)=-\ln\left( 1-2\left(1-\frac{2}{N}\right) t\right) ,\,t\in[0,1/2].
    \]
    One S-type mesh with a similar mesh generating function is the Bakhvalov-S-mesh \cite{RL99} with
    \[
       \phi_{BS}(t)=-\ln\left( 1-2\left(1-\frac{1}{N}\right) t\right) ,\,t\in[0,1/2].
    \]
    The benefit of casting a mesh as S-type mesh is, that using the mesh-characterising function
    \[
       \psi=\exp(-\phi)
    \]
    the convergence rate can be characterised in terms of
    $
       N^{-1}\max|\psi'|
    $
    and optimal meshes are those, where $\max|\psi'|\leq C$ holds independent of $\eps$ and $N$.
    
    Here we have
    \begin{align*}
       \max_{t\in[0,1/2]}|\psi_{BS}'(t)|&=\left|\frac{2}{N}-2\right|\leq 2,&
       \max_{t\in[0,1/2]}|\psi_{eXp}'(t)|&=\left|\frac{4}{N}-2\right|\leq 2.
    \end{align*}
    But $\phi_{eXp}$ does not generate an S-type mesh because
    \[
       \phi_{eXp}(1/2)=\ln \frac{N}{2}
    \]
    while $\phi(1/2)=\ln N$ is needed for an S-type mesh. 
    
    We therefore propose a generalisation of S-type meshes, 
    where the mesh-generating functions $\widehat\phi$ is monotonically increasing and fulfils for a fixed $\alpha>0$
    \[
       \widehat\phi(0)=0,\quad
       \widehat\phi(1/2)=\ln(\alpha N),\quad
       \frac{\max\limits_{t\in[0,1/2]} \widehat\phi'(t)}{N}\leq C.
    \]
    For $\alpha=1$ these are the conditions given in \cite{RL99} for S-type meshes. 
    With
    \[
       \phi_{eXp}'(t)=\frac{2-\frac{4}{N}}{1-2\left(1-\frac{2}{N}\right)t}\leq N-2
    \]
    the given function $\phi_{eXp}$ fulfils these conditions with $\alpha=1/2$ and by using \eqref{eq:Stype} for defining the mesh nodes,
    the exponentially graded eXp-mesh is a modified S-type mesh.    
    A further observation shows for a layer function $E(x)=\e^{-\frac{\beta}{\eps}x}$ at the transition point
    \begin{gather}\label{eq:decay}
       E\left(\frac{\sigma\eps}{\beta}\widehat\phi(1/2)\right)=(\alpha N)^{-\sigma}.
    \end{gather}
    Because only the above properties are needed to prove convergence of FEM on S-type meshes,
    all those proofs apply for the generalised S-type meshes too.
    
    Although \eqref{eq:decay} indicates a better estimation for the standard Bakhvalov-S-mesh (${\alpha=1}$) compared 
    to the eXp-S-mesh ($\alpha=1/2$), numerical results show it vice versa with convergence on the eXp-S-mesh slightly 
    better by a factor, see Section~\ref{sec:numerics}.

  \section{Numerical example}\label{sec:numerics}
    We include only a one-dimensional example. Further examples in two space dimensions showing optimal convergence also compared 
    to Bakhvalov-S-type meshes can be found in \cite{FrLX16,FrLX16_2}. Let us consider the singularly-perturbed
    convection-diffusion problem
    \begin{align*}
      -\eps u''-u'+u&=f,\quad\text{in }(0,1)\\
           u(0)=u(1)&=0,
    \end{align*}
    where the right-hand side is given such that 
    \[
      u(x)=\cos\frac{\pi x}{2}-\frac{\exp(-x/\eps)-\exp(-1/\eps)}{1-\exp(-1/\eps)}
    \]
    is the exact solution for comparing the results. For the numerical example we fix $\eps=10^{-6}$ which is small enough
    to bring out the layer features. Furthermore, $\sigma=p+1$ for polynomial degrees $p=1,\,2,\,3$. 
    Then the discrete spaces for our Galerkin FEM are piecewise polynomials of degree $p$.
    
    We consider the following meshes:
    \begin{enumerate}[1)]
      \item Bakhvalov-S-mesh with $N/2$ elements in the layer region,
      \item eXp-S-mesh according to \eqref{eq:Stype} with $\phi_{eXp}$ and $N/2$ elements in the layer region,
      \item eXp-mesh according to \eqref{eq:eXp},
      \item Bakhvalov-S-mesh with $N/2-1$ elements in the layer region,
      \item eXp-S-mesh with $N/2-1$ elements in the layer region.
    \end{enumerate}
    The reason for including the last two mesh-variations is that here the eXp-mesh and the eXp-S-mesh are identical.
    Also for comparing the quality of the mesh it makes a difference, how many cells are in the layer region.
    
    Table \ref{tab:all}
    \begin{table}[ht]
      \caption{Comparison of the results on all meshes\label{tab:all}}
      \begin{center}   
        \begin{tabular}{rllllllllll}
          $N$  &
          \multicolumn{2}{c}{B-S-mesh} & 
          \multicolumn{2}{c}{eXp-S-mesh} &
          \multicolumn{2}{c}{eXp-mesh} &
          \multicolumn{2}{c}{B-S-mesh$^*$} & 
          \multicolumn{2}{c}{eXp-S-mesh$^*$}\\
          \hline
          \multicolumn{11}{c}{$p=1$}\\
          \hline
           64 & 1.776e-02 &      & 1.749e-02 &      & 1.804e-02 &      & 1.833e-02 &      & 1.804e-02 &     \\
          128 & 8.951e-03 & 0.99 & 8.881e-03 & 0.98 & 9.021e-03 & 1.00 & 9.093e-03 & 1.01 & 9.021e-03 & 1.00\\
          256 & 4.493e-03 & 0.99 & 4.475e-03 & 0.99 & 4.511e-03 & 1.00 & 4.528e-03 & 1.01 & 4.511e-03 & 1.00\\
          512 & 2.251e-03 & 1.00 & 2.246e-03 & 0.99 & 2.255e-03 & 1.00 & 2.260e-03 & 1.00 & 2.255e-03 & 1.00\\          
          \hline
          \multicolumn{11}{c}{$p=2$}\\
          \hline
           64 & 3.903e-04 &      & 3.798e-04 &      & 4.041e-04 &      & 4.157e-04 &      & 4.041e-04 &     \\
          128 & 9.903e-05 & 1.98 & 9.765e-05 & 1.96 & 1.007e-04 & 2.00 & 1.022e-04 & 2.02 & 1.007e-04 & 2.00\\
          256 & 2.494e-05 & 1.99 & 2.476e-05 & 1.98 & 2.515e-05 & 2.00 & 2.533e-05 & 2.01 & 2.515e-05 & 2.00\\
          512 & 6.256e-06 & 2.00 & 6.234e-06 & 1.99 & 6.282e-06 & 2.00 & 6.305e-06 & 2.01 & 6.282e-06 & 2.00\\          
          \hline
          \multicolumn{11}{c}{$p=3$}\\
          \hline
           64& 8.304e-06 &      & 7.938e-06 &      & 8.727e-06 &      & 9.134e-06 &      & 8.727e-06 &     \\
          128& 1.062e-06 & 2.97 & 1.038e-06 & 2.93 & 1.088e-06 & 3.00 & 1.114e-06 & 3.04 & 1.088e-06 & 3.00\\
          256& 1.344e-07 & 2.98 & 1.328e-07 & 2.97 & 1.360e-07 & 3.00 & 1.376e-07 & 3.02 & 1.360e-07 & 3.00\\
          512& 1.689e-08 & 2.99 & 1.679e-08 & 2.98 & 1.699e-08 & 3.00 & 1.709e-08 & 3.01 & 1.699e-08 & 3.00\\          
          \hline
        \end{tabular}
      \end{center}
    \end{table}
  shows the comparison of errors for on all meshes. For all of them the optimal 
  convergence rate of $p$ can be observed for the energy norm
  \[
     \tnorm{u}^2:=\eps\norm{\grad u}{0}^2+\norm{u}{0}^2,
  \]
  where $\norm{\cdot}{0}$ is the usual $L^2$-norm.
  Also the errors using the eXp-S-mesh are the smallest, if only by a small factor. 
  Between the first three meshes the errors on the eXp-mesh are the largest, but 
  the number of cells in the layer-region is different to the S-type meshes. 
  The last three meshes all have $N/2-1$ cells in the layer region.
  Now the results for the eXp-mesh and the eXp-S-mesh$^*$, where one cell less was 
  used in the layer region, are the same because these are numerically the same meshes. 
  And again, the results are slightly better than for the Bakhvalov-S-mesh with $N/2-1$ 
  cells in the layer region.
  
  The improvement of the error on eXp-meshes compared to Bakhvalov-S-meshes is even 
  stronger in 2d, see e.g. the results in \cite{FrLX16}.
 
\section*{Conclusion}
  In this short note we generalised the class of S-type meshes and showed the eXp-mesh 
  to be such a modified S-type mesh. Within the framework of S-type meshes we were 
  able to prove the eXp-S-mesh to be an optimal mesh, allowing a convergence behaviour
  of numerical methods to be slightly better than the Bakhvalov-S-mesh.
  
  \bibliographystyle{plain}
  \bibliography{lit}

\end{document}